\documentclass[12pt]{amsart}
\usepackage{amsfonts,amssymb,latexsym,amsmath, amsxtra}
\usepackage[all]{xy}
\usepackage[dvips]{graphics}

\usepackage[OT2,T1]{fontenc}
\DeclareSymbolFont{cyrletters}{OT2}{wncyr}{m}{n}
\DeclareMathSymbol{\Sha}{\mathalpha}{cyrletters}{"58}

\theoremstyle{plain}
\newtheorem{theorem}{Theorem}[section]
\newtheorem{corollary}[theorem]{Corollary}

\newtheorem{lemma}[theorem]{Lemma}

\newtheorem{proposition}[theorem]{Proposition}

\newtheorem*{conjecture*}{Conjecture}

\theoremstyle{definition}
\newtheorem{definition}[theorem]{Definition}
\theoremstyle{remark}

\numberwithin{equation}{section}
\newcommand{\be}{\begin{equation}}
\newcommand{\ee}{\end{equation}}

\let\non\nonumber

\newcommand{\R}{\mathbb R}
\newcommand{\N}{\mathbb N}
\newcommand{\Z}{\mathbb Z}
\newcommand{\C}{\mathbb C}

\newcommand{\bP}{\mathbb{P}}
\newcommand{\CN}{\mathcal{N}}
\newcommand{\leg}[2]{\left( \frac{#1}{#2} \right)}
\newcommand{\CM}{\mathcal{M}}

\newcommand{\half}{\textstyle{\frac{1}{2}}}
\newcommand{\sgn}{\mathrm{sgn}}


\begin{document}

\title[From sheaves on $\mathbb{P}^2$ to a generalization of the Rademacher expansion]{From sheaves on $\mathbb{P}^2$ to a generalization of the Rademacher expansion}
\author{Kathrin Bringmann}
\address{Mathematical Institute\\University of
Cologne\\ Weyertal 86-90 \\ 50931 Cologne \\Germany}
\email{kbringma@math.uni-koeln.de}

\author{Jan Manschot}
\address{Institut de Physique Th\'eorique\\
 CEA Saclay, CNRS-URA 2306\\
91191 Gif sur Yvette, France}
\email{manschot@uni-bonn.de}

\begin{abstract}
Moduli spaces of stable coherent sheaves on a surface are of much interest for
both mathematics and physics. Yoshioka  computed generating
functions of Poincar\'e polynomials of such moduli spaces if the
surface is $\bP^2$ and the rank of the sheaves is 2. Motivated by physical
arguments, this paper investigates the modular properties of these
generating functions.  It is shown that these functions can be written
in terms of the Lerch sum and theta function. Based on this, we prove a conjecture by Vafa and Witten, which expresses
the generating functions of Euler numbers as a mixed mock modular
form. Moreover, we derive an exact formula of Rademacher-type for the Fourier
coefficients of this function. This formula requires a generalization of
the classical Circle Method. This is the first example of an exact formula for the
Fourier coefficients of mixed mock modular forms, which is of independent mathematical interest.
\end{abstract}
\maketitle

\section{Introduction and Statement of Results}

In the past interactions between physics and mathematics have led to  many
interesting results. Motivated by strong-weak coupling duality (or $S$-duality) in
physics, this article considers various generating functions which
appear in the study of moduli spaces of stable coherent sheaves on the
projective plane $\mathbb{P}^2$. We express the generating functions of Poincar\'e polynomials of moduli spaces of rank 2 sheaves in terms of the  Lerch
sum and theta function which we will recall later. Using these expressions, we prove a conjecture by
Vafa and Witten \cite{Vafa:1994tf} for the generating functions of Euler numbers. These functions appear to be related to Ramanujan's mock theta
functions  and therefore transform almost as weakly holomorphic modular forms,
i.e.,  meromorphic modular forms whose poles (if there are any) may only lie in   cusps.
Our second main result is an exact formula  for the Fourier coefficients of these
generating functions that formally resembles  the Rademacher expansion for the
coefficients of weakly holomorphic modular forms.

Moduli spaces of coherent sheaves on a complex surface $S$ receive much
attention (see for example \cite{Huybrechts:1996} for an extensive work
on such moduli spaces). More specifically, one is interested in the moduli space
$\mathcal{M}(r,c_1,c_2)$ of semi-stable sheaves of rank $r$ with first
Chern class $c_1$ and second Chern class $c_2$. We will
consider topological invariants of $\mathcal{M}$, in particular the Poincar\'e polynomial
$p(\mathcal{M},s):=\sum_{i=0}^{2\dim_{\mathbb{C}}\CM}b_i(\mathcal{M})\,s^i$ and the Euler number
$\chi(\mathcal{M}):=p(\mathcal{M},-1)$, where
$b_i(\mathcal{M})$ is the $i$th Betti number: $b_i(\mathcal{M}):=\dim
H_i(\mathcal{M},\mathbb{Z})$. Ellingsrud and Str\o
mme \cite{Ellingsrud:1987} computed the Betti numbers of the moduli
space of sheaves with rank 1 on $\bP^2$ and other ruled surfaces.  G\"ottsche
\cite{Gottsche:1990} derived the generating function for $p(\mathcal{M}(1,0,n),s):=\sum_{i=0}^{4n }b_i(\mathcal{M})s^i$ for rank 1
sheaves on a smooth projective surface, and wrote it as an elegant
product formula (\ref{eq:poincarer1}).

Subsequent work by Yoshioka \cite{Yoshioka:1994, Yoshioka:1995} derived the generating functions of Poincar\'e polynomials for sheaves
of  rank 2 on the projective plane $\mathbb{P}^2$ (Eqs. (\ref{eq:poincarepols1}) and (\ref{eq:poincarepols0})). Recently, the generating functions of
the Euler and Betti numbers for rank 3  have also been computed \cite{Kool:2009, Manschot:2010nc}.
 Closely related are the computations of the Euler numbers of the moduli spaces of vector bundles for  rank 2
  \cite{Klyachko:1991} and for rank 3  \cite{weist:2009}.

These generating functions have also enjoyed much interest in physics,
in particular in the context of strong-weak coupling duality and instanton moduli
spaces. The duality was  first conjectured by Montonen and Olive
\cite{Montonen:1977sn} as a duality of gauge theory. Their
conjecture claims that gauge theory with gauge group $G$ and coupling
constant $g$ has a dual description in terms of the gauge theory
with gauge group $^L G$ (the Langlands dual group) and coupling constant $4\pi/g$. If the
theta angle $\theta$ is included in the analysis, then this $\mathbb{Z}/2$ group is enlarged
to the modular group SL$_2(\mathbb{Z})$ which acts on the complex parameter
$\tau=\frac{\theta}{2\pi}+\frac{4\pi i}{g^2}\in \mathcal{H}$ by linear fractional
transformations. The symmetry of gauge theory under this larger group
is known as $S$-duality.

Vafa and Witten \cite{Vafa:1994tf} have tested $S$-duality for
topologically twisted gauge theory with $\CN=4$
supersymmetry. They showed that for fixed instanton
number the path integral of this theory equals the Euler number of a suitable
compactification of the instanton moduli space, which  turns out to be
the Gieseker-Maruyama compactification of the moduli space of semi-stable
sheaves whose Chern classes are determined by the instanton
data. $S$-duality led  Vafa and Witten \cite{Vafa:1994tf} (see in
particular Section 3 of \cite{Vafa:1994tf}) to the conjecture that the generating
function of the Euler numbers transforms as a (weakly holomorphic) modular
form with a specific weight and multiplier. These
properties were later also understood from the point of
view of M5-branes, see for example \cite{Minahan:1998vr}.

The generating functions (\ref{eq:poincarepols1}) and
(\ref{eq:poincarepols0}) allow a precise test of the conjectured modular
properties for rank 2 and $\bP^2$. In fact, we derive similar modular properties for these
generating functions of Poincar\'e polynomials as for those of Euler
numbers. This is quite remarkable since present discussions in the literature
are limited to the Euler numbers. To make the modular properties
manifest, we express in Proposition \ref{prop:lerch} these functions in terms of automorphic functions, in particular the Lerch
sum (\ref{eq:lerch}). Since the modular properties of the Lerch sum are  well established, thanks to
Zwegers' thesis \cite{Zwegers:2000}, it is straightforward to derive the modular
properties of these generating functions.

Specialization of the Poincar\'e polynomials to the Euler numbers
requires one to take the derivative of the Lerch sum.
Using this relation we prove that the generating function of the Euler
numbers contains the generating function of the Hurwitz class numbers.A connection between the Euler
numbers of the moduli space of vector bundles and class numbers was earlier proposed by Klyachko
\cite{Klyachko:1991}.  To state our result,
let  $H(n)$ be the Hurwitz class number,   i.e., the number of equivalence classes of quadratic forms of discriminant $-n$, where each class $C$ is
counted with multiplicity $1/\text{Aut}(C)$. We note that $H(0)=-\frac{1}{12}\text{ and } H(3)=\frac{1}{3}$.
Moreover, we let (throughout $q:=e^{2 \pi i \tau}$)
\begin{equation} \label{classgen}
h_j(\tau):=\sum_{n=0}^{\infty}H(4n+3j)q^{n+\frac{3j}{4}}, \qquad j \in\{0,1\}.
\end{equation}
In Section 2, we prove:
\begin{proposition}
\label{col:class}
The generating functions of the  Euler numbers $\chi\left(\mathcal{M}(2,c_1,c_2)\right)$ take the form:
\begin{eqnarray}
\label{eq:P20}
q^{-\frac{1}{2}}\sum_{n=1}^\infty
\chi\left(\mathcal{M}(2,-1,n)\right) q^n&=&\frac{3h_1(\tau)}{\eta^6(\tau)},\non\\
q^{-\frac{1}{4}}\sum_{n=2}^\infty
\chi\left(\mathcal{M}(2,0,n)\right)q^n
&=& \frac{3h_0(\tau)}{\eta^6(\tau)}+\frac{1}{4\eta^3(2\tau)},\non
\end{eqnarray}
where  $\eta(\tau):=q^{\frac{1}{24}}\prod_{n=1}^{\infty}\left(1-q^n
\right)$ is Dedekind's eta-function.
\end{proposition}
\noindent This proposition is the conjecture \cite{Vafa:1994tf} mentioned in the abstract. Results of
Refs. \cite{Klyachko:1991, Yoshioka:1994} led Ref. \cite{Vafa:1994tf} to this conjecture,
and it was verified by a comparison of the first coefficients of
$3h_1(\tau)/\eta^6(\tau)$ with Eq. (\ref{eq:poincarepols1}). The good
modular properties of $h_j(\tau)$ after addition of a suitable non-holomorphic
term (see Eq. (\ref{eq:genclassj})) was a strong confirmation of the $S$-duality conjecture.
We refer to Sec. 4.2 of Ref. \cite{Vafa:1994tf} for more details.


In the following we recall in more detail what is known about modularity  of generating functions of class numbers of imaginary quadratic fields.
Recall that the Fourier coefficients $r(n)$  of $\Theta_0^3$,  with $\Theta_0(\tau):=\sum_{n \in \Z} q^{n^2}$,
themselves encode class numbers. To be more precise,   by a result of Gauss, we have that
 \begin{equation} \label{Gauss}
 r(n)= 12 \left( H(4n)-2H(n)\right).
 \end{equation}
 If one wants to study   the full generating function for the Hurwitz class numbers, then one has to move to the world of harmonic
 weak Maass forms \cite{BF} (see also Section \ref{sec:genfunctions}). These
  are generalizations of modular forms in that they satisfy the same modular transformation laws but instead of being meromorphic
   they are annihilated by the weight $k$ hyperbolic Laplacian.
 To be more precise, Zagier \cite{zagier:1975} showed that  the generating function
 $$
 h(\tau)
 :=  \sum_{\begin{subarray}{c} n\ge 0 \\ n\equiv 0, 3\pmod{4} \end{subarray}} H(n)q^n,
 $$
is  a  mock modular form with shadow $\Theta_0(\tau)$, notions which we will recall shortly.
We note that the full generating function for class numbers of real and imaginary quadratic fields
 requires one to consider even more generalized automorphic objects
\cite{DIT}.

Mock modular forms are related to  Ramanujan's
so-called mock theta functions, which he  introduced  in  his last letter to Hardy (see \cite{Ra}, pp. 127-131) in  17 examples.
  Ramanujan stated that these forms have properties which resemble those of theta functions, but are  not modular forms.
   The mock theta functions occur on the one hand in  a vast variety of papers (see for example  \cite{An, AH,Ch,Hi,Wa} just to mention a few), but
   were on the other hand not well understood for a long time  since they lack real modularity properties.
   The mystery surrounding these functions was finally solved by Zwegers in his famous PhD thesis  \cite{Zwegers:2000}  in which he related the mock theta functions to  harmonic weak Maass forms.
    Placing the mock theta functions into the  world of   harmonic weak Maass forms has many applications: for example the
    first author and Ono proved  an exact formula for the coefficients of one of the mock theta functions \cite{BO1} and explained how to construct  an infinite family of mock theta functions related to Dyson's rank statistic on partitions \cite{BO2}.
   Further applications are for example   a relation between Hurwitz class numbers and overpartition rank differences \cite{BL},
    and a duality relating the coefficients of mock theta functions to coefficients of weakly holomorphic modular forms \cite{FO,Zw}.
    Part of the difficulty of really understanding the mock theta functions was grounded in the fact that these functions have a certain hidden companion,
    which Zagier calls the shadow of the mock theta function, and without which the mock theta functions are not fully understood.
    These shadows  may be obtained from the associated harmonic Maass form by applying the differential operator
    $\xi_{2-k}:=2iy^{2-k} \overline{\frac{\partial}{\partial
        \overline{\tau}}}$ (with $k=1/2$ and $y:=$ Im
        $(\tau)$) and turn out to be unary theta
    functions. Mock modular forms are then  generalizations of mock theta functions in that the associated shadow does not necessarily have
     to be a unary theta function but may be a general (weakly
     holomorphic) modular form. Mixed mock modular forms are functions
     which lie in the tensor space of mock modular forms and modular
     forms.

The functions
\begin{equation} \label{mixedmock}
f_j(\tau):=\frac{h_j(\tau)}{\eta^6(\tau)}=\sum_{n=0}^\infty
\alpha_j(n) q^{n-\frac{j+1}{4}}, \qquad j \in\{0,1\},
\end{equation}
of Proposition \ref{col:class}
are examples of such forms.
In this paper we prove an exact formula for $\alpha_j(n)$, which is the first exact formula for coefficients of mixed mock modular forms.
This result provides an exact formula for
$\chi(\mathcal{M}(2,c_1,c_2))$, which is clearly of
interest for both mathematics and physics.
The proof   requires a generalization of
the Hardy-Ramanujan Circle Method due to the first author and Mahlburg \cite{BM} which may be applied to mixed mock modular forms.
Let us next place this  result in its mathematical context. As usual we denote by $p(n)$ the number of partitions of an integer $n$.
Recall that Hardy and Ramanujan \cite{HR1,HR2}, in work  which gave birth to the Circle Method,
derived their famous asymptotic formula for the partition function $p(n)$,
\begin{equation*}
p(n)\sim \frac{1}{4n\sqrt{3}}\cdot e^{\pi\sqrt{2n/3}} \qquad(n \to \infty).
\end{equation*}
Rademacher  \cite{Rad} then subsequently  proved the following
exact formula
\begin{equation*}
p(n)= 2 \pi (24n-1)^{-\frac{3}{4}} \sum_{k =1}^{\infty}
\frac{A_k(n)}{k}\cdot  I_{\frac{3}{2}}\left( \frac{\pi
\sqrt{24n-1}}{6k}\right).
\end{equation*}
Here $I_{\ell}(x)$
is the $I$-Bessel function of order $\ell$, and $A_k(n)$
is the Kloosterman sum
\begin{displaymath}
A_k(n):=\frac{1}{2} \sqrt{\frac{k}{12}} \sum_{\substack{x \pmod {24k}\\
x^2 \equiv -24n+1 \pmod{24k}}}  \chi_{12}(x) \cdot
e\left(\frac{x}{12k}\right),
\end{displaymath}
where $e(\alpha):=e^{2\pi i \alpha}$ and $\chi_{12}(x):=\leg{12}{x}$.
An important tool used  to prove the asymptotic and exact formulas for $p(n)$ is the fact that
\begin{displaymath}
P(\tau):=\sum_{n=0}^{\infty}p(n)q^{n-\frac{1}{24}}=
\frac{1}{\eta(\tau)}
\end{displaymath}
is a weight $-1/2$
modular form.   Rademacher and Zuckerman
\cite{RZ, Zu1, Zu2} subsequently showed   exact formulas for the
coefficients of generic weakly holomorphic modular forms of negative
weight.

The situation is more complicated if one turns to non-modular objects. Let us mention
 Ramanujan's mock theta functions and in particular
\begin{equation*}
f(q)=\sum_{n=0}^{\infty}\alpha(n)q^n:=
1+\sum_{n=1}^{\infty}\frac{q^{n^2}}{(1+q)^2(1+q^2)^2\cdots
(1+q^n)^2}.
\end{equation*}
The problem of obtaining an asymptotic formula for $\alpha(n)$
is greatly complicated by the fact that $f$ is not a modular form.
Dragonette \cite{Dr} in her PhD thesis   confirmed a conjecture of Ramanujan
concerning an asymptotic formula for $\alpha(n)$;
subsequently this was improved by  Andrews \cite{An2}.
Infinite families of further asymptotic formulas were recently proven by
the first author \cite{Br}.
Andrews and Dragonette moreover  conjectured an exact formula for $\alpha(n)$ which was then
proved by the first author and Ono  \cite{BO1} using the theory of   Maass Poincar\'e series.
The authors of \cite{BO3}  obtained more generally exact formulas for all coefficients of
mock modular forms of non-positive weight.
From the above description it becomes clear that currently asymptotic/exact formulas for coefficients of
modular or mock modular forms are well understood. The situation is totally different for mixed mock modular forms.
The above mentioned methods cannot be applied as the space of harmonic Maass forms is not
closed under multiplication. The first such example was considered by the first author and Mahlburg \cite{BM}
and is  related to so-called partitions without sequences \cite{An3}, a partition statistic that we do not want to recall for the purpose of this paper.
Developing a generalization of the Circle Method to involve certain non-modular objects, the authors managed to obtain asymptotic expansions for such
 partitions without sequences. So far this is the only example of an asymptotic formula for coefficients of forms in the tensor space. It is of mathematical interest to find further such examples.
We note that due to the more complicated situation the authors of \cite{BM} only obtain an asymptotic and not an exact formula.
In this paper we derive the first example of such an exact formula.

Turning back to an exact formula for the coefficients $\alpha_j(n)$ of $f_j$, we require some more notation.
We let for $k \in \N, g\in \Z$, and $u \in \R$
\begin{equation}\label{fkg}
f_{k, g} (u):=
\begin{cases}
\frac{\pi^2}{\sinh^2\left(\frac{\pi u}{k}-\frac{\pi ig}{2k}\right)}  &  \quad \text{if } g\not \equiv 0 \pmod{2k}, \\
\frac{\pi^2}{\sinh^2\left(\frac{\pi u}{k}\right)}-\frac{k^2}{u^2}    &  \quad \text{if } g\equiv 0 \pmod{2k}.
\end{cases}
\end{equation}
Furthermore we define the Kloosterman sums
\begin{equation*}
K_{j, \ell}(n, m;k):=\sum\limits_{{0\leq h< k}\atop{(h, k)=1}}
\psi_{j \ell}(h,h',k)  e^{-\frac{2\pi i}{k}\left(hn+\frac{h'm}{4}\right)},
\end{equation*}
where $\psi_{j\ell}$ is a multiplier defined in (\ref{psimult}), and
$h'$ is given by the congruence $hh'=-1\,(\mod k)$.  Finally we let
$$
\mathcal{I}_{k,g}(n):=
\int_{-1}^{1} f_{k, g}\left(\frac{u}{2}\right)
I_{\frac{7}{2}}\left(\frac{\pi}{k}\sqrt{\left(4n-(j+1)\right)\left(1-u^2\right)}\right)\left(1-u^2\right)^{\frac{7}{4}}du.
$$
\begin{theorem} \label{th:coefficients}
The coefficients $\alpha_j(n)$ of $f_j$ are given by the following exact formula:
\begin{align*}
& \alpha_j(n)=
-\frac{\pi}{6} \left(4n-(j+1)\right)^{-\frac54} \sum_{k=1}^\infty
\frac{K_{j,0}(n, 0;k)}{k}
I_{\frac{5}{2}}\left(\frac{\pi}{k}\sqrt{4n-(j+1)}\right)\\
& +\frac{1}{\sqrt{2}}   \left(4n-(j+1)\right)^{-\frac32}      \sum_{k=1}^\infty\frac{K_{j,0}(n, 0;k)}{\sqrt{k}}I_3
\left(\frac{\pi}{k}\sqrt{4n-(j+1)}\right) \\
& -\frac{1}{8\pi}\left(4n-(j+1)\right)^{-\frac{7}{4}}\sum_{k=1}^\infty
\sum_{\substack{\substack{\ell \in \{0,1\}\\ -k<g \leq k \\ g \equiv \ell \pmod 2 } }}
\frac{K_{j,\ell}(n, g^2;k)}{k^2} \mathcal{I}_{k,g}(n).
\end{align*}
\end{theorem}
The integrals $\mathcal{I}_{k,g}(n)$ can be estimated using well-known asymptotic formulas for Bessel functions and
Proposition 5.1 of \cite{BM}.
\begin{corollary} \label{mainterms}
The leading asymptotic terms of $\alpha_j(n)$ for $n \to \infty$ are:
$$
\alpha_j(n) =
\left(
\frac{1}{96} n^{-\frac32}  - \frac{1}{32 \pi} n^{-\frac74}  + O \left(n^{-2} \right)
\right)e^{2 \pi \sqrt{n} }.
$$
\end{corollary}
\noindent We want to make two remarks concerning Theorem \ref{th:coefficients}
and Corollary \ref{mainterms}.
\begin{enumerate}
\item Firstly we note that in contrast to mock modular forms, the shadows of the mixed mock modular forms  do contribute to our
leading asymptotic terms. One could  determine  further polynomial lower order main terms.

\item Secondly, we like to mention that the first term in the exact formula are the
coefficients of a negative weight Poincar\'e series as described by
Niebur \cite{niebur:1974}. Numerical experiments by F. Str\"omberg
give strong evidence that generically this term does not converge to an integer.
G. W. Moore and the second author \cite{Manschot:2007ha} considered
negative weight Poincar\'e series (which are essentially a sum over $\Gamma_\infty\backslash$SL$_2(\mathbb{Z})$), because of their interpretation
 in the context of the correspondence between 3-dimensional Anti-de Sitter
space and 2-dimensional conformal field theory \cite{Dijkgraaf:2000fq}.
This led them to consider alternate functions, say $\widetilde{f_j}$, in addition to $f_j$, whose coefficients are given by the first term of Theorem
\ref{th:coefficients}. Interestingly, the theorem shows that $\widetilde{f_j}$
appears naturally as part of $f_j$. Moreover, it is possible to
show that $f_j-\tilde f_j$  can also be written as a sum over
$\Gamma_\infty\backslash$SL$_2(\mathbb{Z})$. We leave a precise
discussion for the future.

\end{enumerate}

\noindent The results of this paper might be relevant for various
applications and current developments, some of which we want to list here:
\begin{enumerate}
\item The appearance of Lerch sums in the generating functions of Poincar\'e
  polynomials is rather intriguing. Besides for $\mathbb{P}^2$, one
  can show, using the results of \cite{Yoshioka:1994, Yoshioka:1995}, that they also appear for rank 2
  sheaves on ruled surfaces. Their appearance is essentially a
  consequence of the contributions of stable bundles to the generating series for a specific polarization. It would be interesting to investigate whether Lerch sums play also a
  role for different surfaces, higher rank sheaves, and related systems like Calabi-Yau black holes.
\item The exact formula for the coefficients of
  $f_{j}$ can be generalized to other mixed mock modular forms.
  Besides the intrinsic mathematical interest, this
  might also prove to be very useful in physics, in particular in discussions on black hole entropy and the
AdS$_3$/CFT$_2$ correspondence. For example the functions $h_{j}/\eta
^{24}$ are known to appear as generating functions of the
degeneracies of $\CN=4$ dyons \cite{Dabholkar:2010}.
\end{enumerate}

The outline of this article is as follows. Section
\ref{sec:genfunctions} reviews briefly the generating functions of
invariants of moduli spaces of stable sheaves, relates those of rank 2
to the Lerch sum and theta function, and  proves the conjecture
by Vafa and Witten. Section \ref{sec:exact} derives the exact formula for the
Fourier  coefficients of $f_{j}$, using
the Hardy-Ramanujan Circle Method.

\section*{Acknowledgements}
The authors wish to thank Atish Dabholkar for helpful conversations
and Frederik Str\"omberg for the numerical experiments mentioned in the remark after Corollary \ref{mainterms}.
The first author was partially supported by NSF grant DMS-0757907 and by the Alfried Krupp prize.
The second author was  partially supported  by the ANR grant BLAN06-3-137168.

\section{Generating functions of topological invariants}
\label{sec:genfunctions}
\setcounter{equation}{0}
Let us start by recalling some of the relevant background. We refer the reader who is unfamiliar with these notions from
  algebraic geometry to consult textbooks like \cite{griffiths, nakahara}.
It is well-known that the Betti numbers of $\mathbb{P}^2$ equal $b_0=b_2=b_4=1$  and
$b_1=b_3=0$, therefore  $\chi(\mathbb{P}^2)=3$. The total Chern class
$c(\mathbb{P}^2)$ of the tangent bundle of $\mathbb{P}^2$ is defined
as:
$$
c\left(\mathbb{P}^2\right):=1+c_1\left(\mathbb{P}^2\right)+c_2\left(\mathbb{P}^2\right)=(1+J)^3,
$$
where $J$ is the hyperplane class and $c_i(\mathbb{P}^2)\in
H^{2i}(\mathbb{P}^2,\mathbb{Z})$. Since these cohomology groups are 1-dimensional, we
will also denote the integrated forms $\int c_i(\mathbb{P}^2)$ by
$c_{i}(\mathbb{P}^2)$, thus $c_1(\mathbb{P}^2)=c_2(\mathbb{P}^2)=3$.

Chern classes $c_i(E)$ are defined for any sheaf $E$ on
$\mathbb{P}^2$, and play a central role in the classification of
sheaves. If no confusion can arise, the Chern classes $c_i(E)$ are in the following
abbreviated by  $c_i$. The complex dimension of the moduli space of stable
sheaves may be written in terms of these Chern classes as
\be
\label{eq:dimmod}
\dim_\mathbb{C}\left( \CM(r,c_1,c_2)\right)=2rc_2-(r-1)c_1^2-r^2+1,
\ee
where $r$ is the rank of the sheaf.
The moduli space of stable sheaves generically depends on the choice
of an ample line bundle over the surface. However, since $b_2(\bP^2)=1$, stability does not
depend on this choice.

%

%



We are interested in the generating functions of the Poincar\'e
polynomials and Euler numbers of moduli spaces $\CM(r,c_1,c_2)$ as
functions of $c_1$ and $c_2$.
Twisting  a sheaf by a line bundle $E\otimes \mathcal{O}(k)$ gives an
isomorphism between the moduli spaces  $\mathcal{M}(r,c_1,c_2)$ and
$\mathcal{M}(c_1+rk,c_2+ (r-1)kc_1 + \frac{1}{2}r(r-1)k^2)$. It is therefore sufficient
to only consider  $c_1 \pmod r$.

The generating function of the Poincar\'e polynomials of
$\CM(1,0,c_2)$ for any surface $S$ is given by \cite{Gottsche:1990}
\be
\label{eq:poincarer1}
\sum_{n\geq 0} p\left(\mathcal{M}(1,0,n),s\right)t^n=\prod_{m\geq
  1}\frac{\left(1+s^{2m-1}t^m\right)^{b_1(S)}\left(1+s^{2m+1}t^m\right)^{b_1(S)}}{\left(1-s^{2(m-1)}t^m\right) ^{b_0(S)}\left(1-s^{2m}t^m\right)^{b_2(S)} (1-s^{2(m+1)}t^m)^{b_0(S)}}.
\ee
To exhibit the modular properties for $S=\mathbb{P}^2$, we write it in terms of the Jacobi theta function
$\theta_1(z;\tau)$, which has the following sum and product expansion ($w=e^{2 \pi iz}$):
\be
\label{eq:theta1}
\theta_1(z;\tau):=i\sum_{r\in \mathbb{Z}+\frac{1}{2}} (-1)^{r-\frac{1}{2}} q^{\frac{r^2}{2}}w^{r}
=i q^{\frac{1}{8}} \left(w^\frac{1}{2}
-w^{-\frac{1}{2}}\right)\prod_{n=1}^{\infty}\left(1-q^n\right)\left(1-wq^n\right)\left(1-w^{-1}q^n\right).
\ee
This theta function   transforms under the generators $S:=\left(\begin{smallmatrix}0&-1\\1&0\end{smallmatrix}\right)$ and $T:=\left(\begin{smallmatrix}1&1\\0&1\end{smallmatrix}\right)$ of SL$_2(\mathbb{Z})$
  as:
\begin{eqnarray*}
 \theta_1\left(\frac{z}{\tau} ; -\frac{1}{\tau} \right)&=& -i\sqrt{-i\tau} \exp\left(\frac{\pi i z^2}{\tau}\right)\,\theta_1(z;\tau),\\
 \theta_1(z;\tau+1)&=&\exp\left(\frac{\pi i}{4} \right)\,\theta_1(z;\tau).
\end{eqnarray*}
Moreover $\theta_1$ has simple zeros at the points $z=n \tau+m$ with $n,m\in \Z$.

With the substitutions $q=s^2t=\exp(2\pi i \tau)$ and
$w=s^2=\exp(2\pi i z)$, equation (\ref{eq:poincarer1}) becomes
\be
\label{eq:poin1theta}
q^{-\frac{1}{8}}\sum_{n\geq 0} p\left(\mathcal{M}(1,0,n),w^\frac{1}{2}\right) \left(qw^{-1}\right)^n=\frac{i\left(w^{\frac{1}{2}}-w^{-\frac{1}{2}}\right)}{\theta_1(z;\tau)}.
\ee
The Betti numbers can be obtained by first expanding equation
(\ref{eq:poin1theta}) in $q$ for $q\approx 0$, and then in $w$ for
$w\approx 0$. One easily sees that  (\ref{eq:poin1theta}) has no poles
for $z\in\mathbb{Z}$, but does have simple poles for $z=m\tau+n$, with
$(m,n)\in \mathbb{Z}^2,\,m\neq 0$. The Fourier coefficients of (\ref{eq:poin1theta}) depend
therefore on the choice of contour to extract the Fourier  coefficients.
The physical origin of these poles is however unclear, but  might be related
to the fact that the Poincar\'e polynomial is not a supersymmetric
index like the Euler number.

The above substitutions  for $s$ and $t$ are not arbitrary but
compatible with the Lefshetz $sl(2)$-action on
the moduli space. If $J_3$ is identified with the Cartan element
$sl(2)$, then the action of $J_3$ on an harmonic form on the moduli
space is  given by \cite{Diaconescu:2007bf, griffiths}
$$
J_3\,\omega=\half\left(\deg \omega - \dim \mathcal{M}\right)\,\omega.
$$
The eigenvalue of $J_3$ is the exponent of $w$ in the expansion.

The Euler characteristics are obtained by setting $s=-1$: $p(\mathcal{M}(1,0,n),-1)=\chi(\mathcal{M}(1,0,n))$. Then the generating function becomes:
\be
\label{eq:r1}
f_{1,0}(\tau):=q^{-\frac{1}{8}}\sum_{n=0}^\infty \chi \left(\mathcal{M}(1,0,n)\right)\, q^n=
\frac{1}{\eta^3(\tau)},
\ee
thus  a modular form of weight $-3/2$. For $r\geq 1$, one can more generally define the functions
\be
\label{eq:vecvalued}
f_{r,c_1}(\tau):=\sum_{c_2 \geq \frac{r-1}{2r}c_1^2}\chi
\left(\mathcal{M}(r,c_1,c_2)\right)\,q^{r\Delta - r\chi(\bP^2)/24},
\ee
with $\Delta$ the discriminant of $E$:
$\Delta:=\textstyle{\frac{1}{r}\left(c_2-\frac{r-1}{2r}c_1^2
  \right)}$ and $0 \leq c_1\leq r-1$.
  These functions are expected to exhibit transformation
properties of a vector-valued modular form of length $r$ and weight
$-\chi(\bP^2)/2$. The modular properties  of this function are most straightforwardly
derived from the point of view of multiple M5-branes wrapping
$\bP^2\otimes T^2$. This leads to a generating function which also
sums over all $c_1$ \cite{Minahan:1998vr, Manschot:2008zb}:
\be
\label{eq:ellgenus}
\mathcal{Z}_r(\rho;\tau):=\sum_{c_1, c_2\in \mathbb{Z}}
\chi(\CM(r,c_1,c_2)) \bar q^{r\left(\Delta-\frac{\chi(\bP^2)}{24}\right)}
q^{\frac{1}{2r}\left(c_1+\frac{rc_1\left(\bP^2\right)}{2}\right)^2} (-\xi)^{c_1+\frac{rc_1\left(\bP^2\right)}{2}}
\ee
with $\xi:=e^{2\pi i\rho}$. Physical arguments
suggest that this function transforms under SL$_2(\mathbb{Z})$ like a
Jacobi form which is non-holomorphic in $\tau$ and has weight
$(\frac{1}{2},-\frac{3}{2})$. Moreover, the isomorphism of moduli
spaces due to twisting by a line bundle implies a decomposition of
$\mathcal{Z}_r(\rho;\tau)$ into theta functions
$\Theta_{r,\mu}(\rho;\tau)$ and  vector-valued modular forms
$f_{r,\mu}(\tau)$:
\be
\label{eq:thetadecomp}
\mathcal{Z}_r(\rho;\tau)=\sum_{\mu \pmod r} \bar f_{r,\mu}(\tau) \Theta_{r,\mu}(\rho;\tau),
\ee
with
$$
\Theta_{r,\mu}(\rho;\tau):=\sum_{n=\mu \pmod r}
q^{\frac{1}{2r}\left(n+\frac{rc_1\left(\bP^2\right)}{2}\right)^2} (-\xi)^{n+\frac{rc_1\left(\bP^2\right)}{2}}.
$$
This  decomposition implies that
 $$
 D_r\left(\mathcal{Z}_r(\rho;\tau)\right)=0
 $$
  with $D_r:=\frac{\partial}{\partial \tau}+\frac{i}{4\pi
    r}\frac{\partial^2}{\partial \rho^2}$. Functions which satisfy
  this condition together with an appropriate transformation law are
  known as skew (weakly) holomorphic Jacobi forms \cite{Skoruppa:1988}. In particular for $r=1$ we have
  that
$$
\mathcal{Z}_1(\rho;\tau)=\frac{\theta_1(\rho;\tau)}{\bar \eta^3(\tau)} .
$$
 Already for rank
2, we will find two refinements of these physical expectations.
As explained in the introduction, the $f_{2,\mu}(\tau)$ appear to be mixed mock
modular forms, such that $D_2\left(\mathcal{ \widehat
    Z}_2(\rho;\tau)\right)\neq 0$. This is in physics called a ``holomorphic anomaly''. The other refinement
concerns the integrality of the Fourier coefficients.

%
%

  Before  returning to the functions of interest for this paper, we
   want to recall the precise definition of harmonic weak Maass forms.
   Here we only require the case of half-integral weight on $\Gamma_0(4)$.
  \begin{definition}
A harmonic weak Maass form of weight  $k \in \frac12 +\Z$ for the group $\Gamma_0(4)$
is a smooth function $f:\mathbb{H} \to \mathbb{C}$
satisfying the following:
\begin{enumerate}
\item
For all $\left(\begin{smallmatrix} a & b \\ c & d \end{smallmatrix} \right)\in \Gamma_0(4)$, we have that
\begin{displaymath}
f\left(\frac{a \tau+b}{c \tau+d}  \right)=
\leg{c}{d}\epsilon_d^{-2k}(c\tau+d)^{k}\ f(\tau) .
\end{displaymath}
Here $\leg{c}{d}$ denotes the Jacobi symbol, $\epsilon_d=1$ for $d\equiv 1 \pmod{4}$
and $\epsilon_d=i$ for $d\equiv 3 \pmod{4}$, and $\sqrt{\tau}$ is the principal branch of the holomorphic square root.
\item We have $\Delta_k f=0$, where $(\tau=x+iy$) the weight $k$ hyperbolic Laplacian $\Delta_k$ is defined as
$$
\Delta_k := -y^2\left( \frac{\partial^2}{\partial x^2} +
\frac{\partial^2}{\partial y^2}\right) + iky\left(
\frac{\partial}{\partial x}+i \frac{\partial}{\partial y}\right).
$$
\item The function $f$
has at most
linear exponential growth at  all the cusps.
\end{enumerate}
\end{definition}
Using this notation,  the function
\begin{equation}
\label{eq:genclass}
\widehat h(\tau):=\sum\limits_{{n=0}\atop{n\equiv 0,3\pmod{4}}}^{\infty}H(n)q^n+\frac{(1+i)}{16\pi}\ \displaystyle\int^{i\infty}_{-\overline{\tau}}\
\frac{\Theta_0(w)}{(\tau+w)^\frac{3}{2}}\ dw
\end{equation}
\noindent
is a harmonic Maass form of weight $\frac{3}{2}$ on $\Gamma_0(4)$ (see \cite{zagier:1975}).
We moreover require the restrictions of $\widehat{h}$  to arithmetic progressions $0, 3 \pmod 4$ (individually). It is not hard to
see that the associated harmonic weak Maass forms are given by
\begin{equation}
\label{eq:genclassj}
 \widehat h_j(\tau):=h_j(\tau)+\frac{(1+i)}{8\pi}\displaystyle\int_{-\overline{\tau}}^{i\infty}
\frac{\Theta_j(w)}{(\tau+w)^{\frac{3}{2}}}\,dw,\qquad j \in \{0,1\},
\end{equation}
where the functions $h_j$ were defined in (\ref{classgen})
and
$\Theta_j(\tau):=\displaystyle\sum_{n\in\Z}q^{\frac{1}{4}(2n+j)^2}.$

Now we continue our discussion on generating functions related to
semi-stable coherent sheaves of rank 2 on $\mathbb{P}^2$.
Yoshioka  \cite{Yoshioka:1994, Yoshioka:1995} computed the
generating functions of the Poincar\'e polynomials. To present his result, define
$$
Z_s\left(\bP^2,t\right):=\frac{1}{(1-t)(1-st)\left(1-s^2t\right)}.
$$
The expression for the generating function for $c_1=-1$ is \cite{Yoshioka:1994}:
\begin{eqnarray}
\label{eq:poincarepols1}
&&\sum_{n=1}^\infty p\left(\mathcal{M}(2,-1,n),s\right)t^n=\frac{\prod_{d\geq 1}Z_{s^2}\left(\bP^2,s^{4d-2}t^d\right)^2}{\left(s^2-1\right)\sum_{n\in
    \mathbb{Z}}s^{2n(2n-1)}t^{n^2}}\\
&&\qquad\times \sum_{b\geq
  0}\left(\frac{s^{2(b+1)(2b+1)}}{1-s^{8(b+1)}t^{2b+1}}-\frac{s^{2b(2b+5)}}{1-s^{8b}t^{2b+1}}
\right)t^{(b+1)^2}, \non
\end{eqnarray}
and similarly for $c_1=0$ \cite{Yoshioka:1995}:
\begin{multline}
\label{eq:poincarepols0}
\sum_{n=2}^\infty p\left(\mathcal{M}(2,0,n),s\right)t^n=\frac{\prod_{d\geq 1}Z_{s^2}\left(\bP^2,s^{4d-2}t^d\right)^2}{\left(1-s^2\right) \sum_{n\in
    \mathbb{Z}}s^{2n(2n+1)}t^{n(n+1)} }\\
\times \left(\sum_{b\geq
  0}-\left(\frac{s^{2(b+1)(2b+3)}}{1-s^{8(b+1)}t^{2b+1}}-\frac{s^{2b(2b+7)}}{1-s^{8b}t^{2b+1}}
\right)t^{b^2+3b+1} 
+\sum_{b\geq 0}\frac{s^{2(b+1)(2b+1)}-s^{2b(2b+1)}}{2s^2}t^{b(b+1)}\right) \\
\qquad +\frac{\prod_{d\geq 1}Z_{s^4}\left(\bP^2,s^{8d-4}t^{2d}\right)}{2s^2\left(1+s^2\right)}.
\end{multline}
Here we have corrected a sign error in Remark 4.6 of
\cite{Yoshioka:1995}.

To simplify the expressions (\ref{eq:poincarepols1}) and
(\ref{eq:poincarepols0}), we make the substitutions $s^4t=q$ and
$s^2=w$, analogous to the substitutions in the rank 1 case.
One finds after a straightforward computation:
\begin{proposition}
\label{prop:lerch}
The generating functions of the Poincar\'e polynomials $p\left(\mathcal{M}(2,c_1,c_2),s\right)$ take the form:
\begin{eqnarray}
\label{eq:P2-1}&&
q^{-\frac{1}{2}}\sum_{n=1}^\infty
p\left(\mathcal{M}(2,-1,n),w^\frac{1}{2}\right)\left(qw^{-2}\right)^n=\\
&&\qquad -\frac{(1-w)}{w^\frac{5}{2}\,\theta^2_1(z;\tau)}\,\mu\left(2z-\tau,\half-\tau-z;2\tau\right),\non
\\
\label{eq:P20}&&
q^{-\frac{1}{4}}\sum_{n=2}^\infty
p\left(\mathcal{M}(2,0,n),w^\frac{1}{2}\right)\left(qw^{-2}\right)^n
=\\
&&\qquad \frac{(1-w)}{w^2\,\theta_1^2(z;\tau)}\,\left(\frac{1}{2}-q^{-\frac{1}{4}}w^{\frac{3}{2}}\,
\mu\left(2z-\tau,\half-z;2\tau\right)\right)-\frac{i(1-w)}{2w^2\,\theta_1(2z;2\tau)}\non,
\end{eqnarray}
where $\mu(u,v;\tau)$  is the Lerch sum defined by
\be
\label{eq:lerch}
\mu(u,v;\tau):=\frac{e^{\pi i u}}{\theta_1(v;\tau)}\sum_{n\in
 \mathbb{Z}}\frac{(-1)^ne^{\pi i(n^2+n)\tau+2\pi i nv}}{1-e^{2\pi i n
   \tau+2\pi i u}},
\ee
with $u,v\in \mathbb{C}$.
\end{proposition}

Moreover, we define
\begin{eqnarray}
\label{eq:f21}
f_{2,1}(z;\tau)&:=&\frac{\left(1-w\right)\, q^{-\frac{1}{4}}}{w\, \theta_1^2(z;\tau) \, \theta_3(z;2\tau)}
\sum_{n\in \mathbb{Z}}\frac{q^{n^2}w^{-n}}{1-q^{2n-1}w^2}, \non \\
\label{eq:f20}
f_{2,0}(z;\tau)&:=&\frac{(1-w)}{w^2\,\theta_1^2(z;\tau)}
\left(\frac{1}{2}+ \frac{q^{-\frac{3}{4}}w^\frac{5}{2}}{\theta_2(z;2\tau)}
\sum_{n\in\mathbb{Z}}\frac{q^{n^2+n}w^{-n}}{1-q^{2n-1}w^2} \right), \non
\end{eqnarray}
with
$$
\theta_2(z;\tau):=\sum_{n\in \mathbb{Z}+\frac{1}{2}}q^{\frac{n^2}{2}}w^n,
\qquad \qquad\theta_3(z;\tau):=\sum_{n\in \mathbb{Z}}  q^{\frac{n^2}{2}}w^{n}.
$$
and
\begin{eqnarray*}
&& g_{1}(z;\tau):=\frac{q^{-\frac{1}{4}}w^{\frac{3}{2}}}{\theta_3(z;2\tau)}
\sum_{n\in \mathbb{Z}}\frac{q^{n^2}w^{-n}}{1-q^{2n-1}w^2}
=-\mu\left(2z-\tau,\half-\tau-z;2\tau\right),\\
&& g_0(z;\tau):=\frac{1}{2}+\frac{q^{-\frac{3}{4}}w^\frac{5}{2}}{\theta_2(z;2\tau)}
\sum_{n\in \mathbb{Z}}\frac{q^{n^2+n}w^{-n}}{1-q^{2n-1}w^2}=\frac{1}{2}-q^{-\frac{1}{4}}w^{\frac{3}{2}}\,
\mu\left(2z-\tau,\half-z;2\tau\right),
\end{eqnarray*}
Similarly to the case of $r=1$  these functions have poles for $z\in
m\tau+n$ with $(m,n)\in \mathbb{Z}^2,\,m\neq 0$.

Since we have now  explicit expressions for the generating functions
of Poincar\'e polynomials for $r=1,2$ at our disposal, it is particularly interesting to investigate their modular
properties. For rank 1, equations (\ref{eq:poin1theta}) and
(\ref{eq:r1}) show that the generating function of Euler
numbers is indeed a weakly holomorphic modular form, whereas the generating function
 for Poincar\'e polynomials transforms as a Jacobi form of weight
$-\frac{1}{2}$ and index $-\frac{1}{2}$ (up to  the
prefactor $w^{\frac{1}{2}}-w^{-\frac{1}{2}}$).

To make the modular properties of (\ref{eq:P2-1}) and
(\ref{eq:P20}) more manifest, we recall some
results of Zwegers' thesis  \cite{Zwegers:2000}.
The Lerch sum (\ref{eq:lerch}) does not transform as a Jacobi form under
SL$_2(\mathbb{Z})$. However, the completed function \cite{Zwegers:2000}
$$
\widehat{\mu}(u,v;\tau):=\mu(u,v;\tau)+\frac{i}{2}R(u-v;\tau)
$$
 transforms as a
multi-variable Jacobi form of  weight $\half$. Here the function
$R(u;\tau)$ is defined by:
$$
R(u;\tau):=\sum_{r\in\mathbb{Z}+\frac{1}{2}} \left(  \sgn(r)
-E\left( (r+a)\sqrt{2 y}\right)\right)
(-1)^{r-\frac{1}{2}}e^{-\pi i r^2\tau-2\pi i r u},
$$
with $a:=\text{Im}(u)/y$, and
$$
E(z):=2 \int_{0}^{z}  e^{-\pi u^2}\,du.
$$
To be more precise, we have that:
\begin{enumerate}
\item For $k,l,m,n\in \mathbb{Z}$, we have that:
$$
\widehat{\mu}(u+k\tau+l,v+m\tau+n;\tau)=(-1)^{k+l+m+n}\,e^{\pi i (k-m)^2\tau+2\pi i (k-m)(u-v)}\widehat{\mu}(u,v;\tau).
$$
\item
For $\gamma=\left( \begin{smallmatrix}a & b \\ c & d \end{smallmatrix}\right) \in \text{SL}_2(\Z)$, we have that:
$$
\widehat{\mu}\left(\frac{u}{c\tau+d},\frac{v}{c\tau+d};\frac{a\tau+b}{c\tau+d}\right)
=v(\gamma)^{-3}(c\tau+d)^\frac{1}{2}\,e^{-\frac{\pi i c(u-v)^2}{c\tau+d}}\,\widehat{\mu}(u,v;\tau),
    $$  with
  $v(\gamma):=\eta\left(\frac{a \tau +b}{c \tau +d} \right)/\left((c\tau+d)^\frac{1}{2}\eta(\tau)\right)$.
\end{enumerate}
Moreover we require the following identity, which allows us to shift parameters in the function $\mu$ ($z \in \mathbb{C}$):
\begin{equation} \label{parshift}
\mu(u+z, v+z;\tau)-\mu(u, v;\tau)=\frac{i\eta^3(\tau)\theta_1(u+v+z;\tau)\theta_1(z;\tau)}{\theta_1(u;\tau)\theta_1(v;\tau)\theta_1(u+z;\tau)\theta_1(v+z;\tau)}.
\end{equation}
We now turn back to the functions $g_j(z;\tau)$ and define their completions as:
\begin{equation*}
\widehat{g}_{j}(z;\tau):=g_j(z;\tau)+\frac{1}{2}R_j(z;\tau),
\end{equation*}
with
\begin{equation*}
R_j(z;\tau):=\sum_{n\in\mathbb{Z}+\frac{j}{2}}\left( \sgn(n)-E\left((2n+3a)\sqrt{y}\right)\right)\,q^{-n^2}w^{-3n}
.
\end{equation*}
We note that
\begin{eqnarray*}
R_1(z;\tau)&=& -iR\left( 3z-\frac12; 2 \tau\right),\\
R_0(z;\tau)&=& -1-iq^{-\frac14} w^{\frac32} R\left( 3z-\tau-\frac12; 2 \tau\right).
\end{eqnarray*}
 Using the above stated transformation properties of $\widehat{\mu}$, one can show that the functions
 $\widehat{g}_{j}$ are invariant under $T^4$, and transform under $S^{-1}T^{-4}S\in\Gamma_0(4)$ as
\begin{equation*}
\widehat{g}_j\left(\frac{z}{4\tau+1};\frac{\tau}{4\tau+1}\right)= (4\tau+1)^\frac{1}{2}\exp\left(2\pi i \frac{-9z^2}{4\tau+1}\right)\widehat{g}_j(z;\tau).
\end{equation*}
Moreover, one can prove with some more work that the function $\widehat g_{j}(z;\tau)$ may be viewed as components of a function that transforms
like a vector valued modular form for  SL$_2(\mathbb{Z})$ of  weight $\frac{1}{2}$ and with
the same multipliers as the $\widehat h_{j}(\tau)$ defined in equation
(\ref{eq:genclassj}).

We observe that if $f_{2,j}(z;\tau)$ would be completed to
$\widehat f_{2,j}(z;\tau)$ by changing $g_{j}(z;\tau)$ to $\widehat
g_{j}(z;\tau)$, they would transform as  Jacobi forms for
$\Gamma_0(4)$ of  weight $-\frac{1}{2}$ and index
$-\frac{13}{4}$, if we ignore the prefactors $(1-w)/w^{2-j}$. The
non-holomorphic parts of $\widehat f_{2,j}(z;\tau)$ might appear
naturally in physics, but precisely how is unknown. Since the functions
$\widehat g_j(z;\tau)$ transform as a modular vector, the function
$ \mathcal{\widehat Z}_2(z,\rho;\tau)=\sum_{j=0,1} \overline{ \widehat f_{2,j}(z;\tau)} \Theta_{2,j}(\rho;\tau)$
transforms with weight $(\frac{1}{2},-\frac{1}{2})$ under
SL$_2(\mathbb{Z})$ (ignoring the prefactors), which can be
understood from physics.

%
%

As expected from physical arguments, the modular
properties improve if one takes the limit
$w^\frac{1}{2}\to-1$. One can derive
straightforwardly that in this case the last term in equation (\ref{eq:P20})
is equal to $\frac{1}{4}\eta^{-3}(2\tau)$, which is a modular form of
$\Gamma_0(4)$ with a non-trivial multiplier. More interesting is that the limit translates to taking
the derivative of the Lerch sums in Eqs. (\ref{eq:P2-1}) and (\ref{eq:P20}).
Proposition \ref{col:class} gives for these terms
$f_{2,j}(\tau):= f_{2,j}(0;\tau)=3h_{j}(\tau)/\eta^6(\tau)$.

Before proving the Proposition \ref{col:class}, we would like to make a couple of  remarks concerning
$f_{2,j}(\tau)$, and $\mathcal{Z}_2(\rho;\tau)$ defined by equation
(\ref{eq:ellgenus}). The completions $\widehat f_{2,j}(\tau)$ can be
obtained from $\widehat f_{2,j}(z;\tau)$, by computing the coefficient
of $z^1$ in the Taylor expansion of $\widehat g_j(z;\tau)$. Due to the
non-holomorphic term, $D_2 \left(\mathcal{\widehat
    Z}_2(\rho;\tau)\right)\neq 0$. One finds
\begin{equation*}
D_2 \left(\mathcal{\widehat Z}_2(\rho;\tau)\right)
=\frac{-3i}{16\pi   y^{3/2}}\frac{\theta_1^2(\rho;\tau)}{\overline{\eta}^6(\tau)},
\end{equation*}
which is proportional to $\mathcal{Z}_1^2(\rho;\tau)$. The authors of \cite{Minahan:1998vr} conjecture  that such an anomaly appears
 generically for $r\geq 2$. We did not find such a factorization in
 the case of  Poincar\'e polynomials, that is to say for $D_2 \left(\mathcal{\widehat
     Z}_2(z,\rho;\tau)\right)$.


The generating functions of the Euler numbers, $\widehat f_{2,1}(\tau)$ and $\widehat
f_{2,0}(\tau)+\frac{1}{4}\eta^{-3}(2\tau)$, do not combine to a
vector-valued modular form because of the term
$\frac{1}{4}\eta^{-3}(2\tau)$. This term
disappears, if we consider the generating
functions of the rational invariants
$$
\overline{\chi}(\Gamma):=\sum_{m\geq 1,\, m|\Gamma}
(-1)^{\mathrm{dim}_\mathbb{C}(\mathcal{M}(\Gamma/m))}\chi(\Gamma/m)/m^2,
$$
where $\Gamma$ represents the data of the sheaf
$(r,c_1,c_2)$. Refs. \cite{Manschot:2010xp,  Vafa:1994tf} give also evidence that the generating function of the rational invariants $\overline{\chi}(\Gamma)$ have better modular properties, than the ones for the
integer invariants $\chi(\Gamma)$. This is the second refinement, alluded to below equation (\ref{eq:thetadecomp}).

On the other hand, the coefficients of
$f_{2,0}(\tau)+\frac{1}{4}\eta^{-3}(2\tau)$
are required to be integers. This can easily be seen
from the arithmetic properties of the functions. To see this, multiply the
function by $\eta(\tau)^6$, which gives
$3h_0(\tau)+\frac{1}{4}\Theta_0^3(\tau+\half)$. Integrality of the
coefficients of this function is manifest, due to the
properties of the class numbers $H(n)$ and $\Theta_0^3(\tau+\half)$.

\subsection*{Proof of Proposition \ref{col:class}}
One could prove the proposition straightforwardly by verifying 1) that the
shadows of $f_{2,j}(\tau)$, as obtained from $f_{2,j}(z;\tau)$, coincide with those of $3h_{j}(\tau)/\eta^6(\tau)$,
and 2) that a specific number (related to the dimension of the space of
associated modular forms) of coefficients agree. Since this is rather
technical, we choose to prove the proposition by relating it to known
expressions in the literature.

We start with the identity for $f_{2, 1}$. The limit $z \to 0$
of $f_{2,1}(z;\tau)$ is finite and leads to differentation of the
Lerch sum:
\begin{equation} \label{f21}
f_{2, 1}(\tau)=-\frac{1}{\eta^6(\tau)}\frac{d}{dw}
\left[\mu\left(2z-\tau, -z-\tau+\frac{1}{2}; 2\tau\right)\right]_{w=1}.
\end{equation}
Using (\ref{parshift}) yields that
\begin{align*}
\mu\left(2z-\tau, -z-\tau+\frac{1}{2}; 2\tau\right)
&=\mu\left(-\tau, -3z-\tau+\frac{1}{2}; 2\tau\right)\\
& \quad+\frac{i\eta^3(2\tau)\theta_1\left(-z-2\tau+\frac{1}{2};2 \tau\right)\theta_1(2z;2\tau)}
{\theta_1(2z-\tau;2 \tau)\theta_1\left(-z-\tau+\frac{1}{2};2 \tau\right)\theta_1( -\tau; 2 \tau)\theta_1\left( -3z-\tau+\frac{1}{2}; 2 \tau\right)}.
\end{align*}
One can prove that the   second summand contributes $\frac{1}{2}\Theta_1^3(\tau)/\eta^6(\tau)$ to (\ref{f21}). Moreover, the contribution from the
first summand is given by
\[
-\frac{3}{2}\frac{q^{-\frac14}}{\eta^6(\tau)\,\Theta_0(\tau)}\sum_{n\in\Z}\frac{(2n-1)q^{n^2}}{1-q^{2n-1}}.
\]
Using work of Kronecker \cite{Kr1}, Mordell \cite{Mo1}, and Watson \cite{Wa2}, one can prove that
\begin{equation} \label{classidentity}
h_1(\tau)=-\frac{1}{2 \Theta_0(\tau)}q^{-\frac{1}{4}}\sum_{n\in\Z}\frac{(2n-1)q^{n^2}}{1-q^{2n-1}}+\frac{1}{6}\Theta_{1}^3(\tau).
\end{equation}
From this the claim may be easily concluded. We first note that by Watson (correcting a typo) we obtain that
\begin{equation} \label{Watson}
\sum_{n=0}^\infty F(4n+3)q^{n+\frac{3}{4}}
= \frac14\Theta_1^3(\tau)-\frac{1}{\vartheta_3(0)}
\sum_{n\in\Z}\frac{\left(n-\frac12\right)q^{\left(n-\frac12\right)^2}}{q^{\frac12-n}-q^{n-\frac12}}
\end{equation}
where $F(n)$ counts the number of uneven equivalence classes of positive definite quadratic forms of discriminant $-n$.
Next one can easily show (for example by using the theory of modular forms) that
\[
\Theta_1^3(\tau)=\sum_{n=0}^\infty r(4n+3) q^{n+\frac34},
\]
where the coefficient $r(n)$ is defined by
\[
\Theta_0^3(\tau)=\sum_{n=0}^\infty r(n) q^n.
\]
Now a direct computation gives (\ref{classidentity}).

We next turn to  $f_{2,0}$. It is not hard to see that
\begin{equation} \label{f2equation}
f_{2, 0}(\tau)=-\frac{1}{\eta^6(\tau)} \frac{d}{dw}\left[q^{-\frac{1}{4}}w^{\frac{3}{2}}\mu\left(2z-\tau, -z+\frac{1}{2}; 2\tau\right)\right]_{w=1}.
\end{equation}
We find using (\ref{parshift})
\begin{multline*}
\mu\left(2z-\tau, -z+\frac{1}{2}; 2\tau\right)=\mu\left(-\frac{1}{2}, -3z+\tau; 2\tau\right)\\
+\frac{i\eta^3(2\tau)\theta_1(-z;2\tau)\theta_1\left(2z-\tau+\frac{1}{2}; 2 \tau\right)}{\theta_1(2z-\tau; 2 \tau)\theta_1\left( -z+\frac{1}{2}; 2 \tau\right)\theta_1\left(-\frac{1}{2}; 2 \tau\right)\theta_1(-3z+\tau; 2 \tau)}.
\end{multline*}
One can show that the contribution of the second summand to (\ref{f2equation}) equals $-\frac14 \Theta_0^3(\tau)/\eta^6(\tau)$.
Moreover one can prove that the first summand gives a contribution of
$$
  \frac{-3}{\eta^6(\tau)\,\Theta_0\left(\tau+ \frac12 \right)}
\sum_{n \in \Z} \frac{n(-1)^n  q^{n^2} }{1+q^{2n}}.
$$
Now the claim easily follows using (\ref{Gauss}) and the identity
\begin{equation}
\label{eq:classidentity}
\sum_{n=0}^{\infty}H(n)q^n= -\frac{1}{2\Theta_0\left(\tau+\frac{1}{2}\right)}\sum_{n\in\Z}\frac{n(-1)^nq^{n^2}}{1+q^{2n}}
- \frac{1}{12} \Theta_0^3(\tau) .
\end{equation}
Indeed, equation (\ref{eq:classidentity}) may for example be concluded by combining Theorem 1.1 and
Corollary 1.6 of \cite{BL1} and inserting the generating function for $\overline{f}$
given in \cite{BL}. $\square$
\vspace{.2cm}

\section{Exact formulas for  $\alpha_j(n)$}
\label{sec:exact}
The introduction motivates the derivation of an exact formula of
the Fourier coefficients of $f_0 (\tau)=\frac{1}{3}f_{2,0}(\tau)$ and $f_1 (\tau)=\frac{1}{3}f_{2,1}(\tau)$. This will be the subject of this
section. We start by providing various useful transformation formulas,
after which we use the Hardy-Ramanujan Circle Method, to derive the exact
formula.

\subsection{Some transformation formulas} \label{TransformationSection}
In this section, we give transformation properties for the class number generating functions $h_0$ and $h_1$.
Throughout, we let $z \in \C$ with Re$(z)>0$, $k>0$, $(h,k)=1$, and $h'$ defined
via the congruence $hh'\equiv -1 \pmod k$.
Moreover, we assume that $4|h'$ if $k$ is odd.
%
We require  the transformation law of the eta-function:
\begin{equation} \label{EtaTrans}
\eta \left( \frac{1}{k}(h+iz)\right)
= e^{ \frac{\pi i }{12k} (h-h') } \cdot
\omega_{h,k}^{-1} \cdot z^{-\frac{1}{2}}  \cdot
\eta \left( \frac{1}{k}\left(h'+\frac{i}{z}\right)\right),
\end{equation}
where $\omega_{h,k}$ is given by
\begin{eqnarray*}
\omega_{h,k} :=
\exp\left(\pi i
 \sum_{\mu \pmod k}  \left( \left( \frac{\mu}{k}\right) \right)   \left( \left(
\frac{h \mu}{k}\right) \right)
  \right),
\end{eqnarray*}
with
\begin{eqnarray*}
((x)):= \left \{
\begin{array}{ll}
x- \lfloor x \rfloor - \frac{1}{2} &\text{if } x \in \R \setminus \Z ,\\
0&\text{if } x \in \Z.
\end{array}
\right.
\end{eqnarray*}
Writing $\Theta_0$ and $\Theta_1$ as eta-quotients
$$
\Theta_0(\tau) = \frac{\eta^5(2 \tau)}{\eta^2(\tau) \eta^2(4 \tau)},
\qquad \qquad
\Theta_1(\tau)=2 \frac{\eta^2(4 \tau)}{\eta(2\tau)}
$$
yields the following  transformation law ($j \in \{0,1\}$):
\begin{equation*}
\Theta_j \left(\frac{1}{k}(h+iz)\right)=\frac{1}{\sqrt{z}}
\sum_{\ell \in \{0,1\}}
\chi_{j \ell}(h, h',k) \Theta_{\ell}\left(\frac{1}{k}\left(h'+\frac{i}{z}\right)\right)
 .
\end{equation*}
Here the multipliers $\chi_{j \ell}$ are defined as follows:
\begin{align*}
\chi_{00}(h, h', k)&:=
\begin{cases}
\frac{\omega_{h,k}^2\omega_{h,\frac{k}{4}}^2}{\omega_{h,\frac{k}{2}}^5}  &\text{
if } 4|k,\\
\frac{1}{\sqrt{2}}   \frac{\omega_{h,k}^2\omega_{4h,k}^2}{\omega_{2h,k}^5}
&\text{ if }  2\nmid k,\\
0 &\text{ if } 2\|k,
\end{cases}  \qquad
\chi_{01}(h, h', k) :=
\begin{cases}
0 &\text{ if } 4|k,\\
 \frac{1}{\sqrt{2}}   \frac{\omega_{h,k}^2\omega_{4h,k}^2}{\omega_{2h,k}^5}
&\text{ if }  2\nmid k,\\
\frac{\omega_{h,k}^2\omega_{2h,\frac{k}{2}}^2}{\omega_{h,\frac{k}{2}}^5}e^{-\frac{\pi
    i h'}{2k}}
 &\text{ if }  2\|k,
\end{cases} \\
\chi_{10}(h, h', k)&:=
\begin{cases}
0                                                                                                                                                                                                                                                                                                                                 &\text{if } 4|k,\\
\frac{1}{\sqrt{2}}         \frac{\omega_{2h,k}}{\omega_{4h,k}^2}e^{\frac{\pi
ih}{2k}}                &\text{if } 2\nmid k,\\
                         \frac{\omega_{h,\frac{k}{2}}}{\omega_{2h,\frac{k}{2}}^2}e^{\frac{\pi
ih}{2k}}                                                                                        &\text{if } 2\|k,
\end{cases} \qquad
\chi_{11}(h, h', k):=
\begin{cases}
 \frac{\omega_{h,\frac{k}{2}}}{\omega_{h,\frac{k}{4}}^2}e^{\frac{\pi
i}{2k}(h-h')}                                                        &\text{if } 4|k,\\
                -\frac{1}{\sqrt{2}}         \frac{\omega_{2h,k}}{\omega_{4h,k}^2}e^{\frac{\pi
ih}{2k}}                                &\text{if } 2\nmid k,\\
0                                                                                                                                                                                                                                                                                                                                         & \text{if } 2\|k.
\end{cases}
\end{align*}
We next use the   well-known behavior of $\widehat h_j$ under inversion and
translation
\begin{align*}
\widehat{h}_0(\tau +1) = \widehat{h}_0(\tau)&\qquad\qquad \widehat{h}_1(\tau+1)=-i \widehat{h}_1(\tau)\\
\widehat{h}_0\left(-\frac{1}{\tau}\right) = \tau^{\frac{3}{2}}\frac{(1+i)}{2}\left(\widehat{h}_0(\tau)+\widehat{h}_1(\tau)\right)&\qquad\qquad
\widehat{h}_1\left(-\frac{1}{\tau}\right) = \tau^{\frac{3}{2}}\frac{(1+i)}{2}\left(\widehat{h}_0(\tau)-\widehat{h}_1(\tau)\right).
\end{align*}
This gives that $\widehat{h}_j$ has  multiplier dual to the one of $\Theta_j$. To be more precise we have
$$
\widehat h_j \left( \frac{1}{k} (h+iz)\right)
= -z^{-\frac32}
\sum_{\ell \in \{0,1\}}
\overline{\chi_{j\ell}(h,h',k)} \,  \widehat h_{\ell} \left( \frac{1}{k}
\left(h'+\frac{i}{z}\right)\right).
$$
From this, a straightforward calculation shows that
\begin{multline} \label{mockT}
h_j\left( \frac{1}{k} (h+iz)\right)
= -z^{-\frac32}
\sum_{\ell \in \{0,1\}}   \overline{\chi_{j \ell}(h,h',k)}   h_{\ell}\left( \frac{1}{k}
\left(h'+\frac{i}{z}\right)\right) \\
 -\frac{1}{4 \sqrt{2}\pi}   z^{-\frac32} \int_{0}^{\infty} %
\frac{\sum_{\ell \in \{0,1\}}  \overline{\chi_{j \ell}(h,h',k)} \Theta_{\ell} \left( it-\frac{h'}{k} \right)}{\left(
t+\frac{1}{kz}\right)^{\frac32}} dt.
 \end{multline}
Defining
\[
\mathcal{I}_j(x) :=
\displaystyle\int_0^{\infty}\frac{\Theta_j(iw-\frac{h'}{k})}{(w+x)^{\frac{3}{2}}}\,dw
\]
we may  rewrite (\ref{mockT}) as
\begin{align*}
 h_j\left(\frac{1}{k}(h+iz)\right)=&
-z^{-\frac{3}{2}}
\sum_{\ell \in \{0,1\}}
\overline{\chi_{j\ell}(h,h', k)}h_{\ell}\left(\frac{1}{k}\left(h'+\frac{i}{z}\right)
\right)
\\
& -\frac{1}{4 \sqrt{2}\pi}
z^{-\frac{3}{2}} \sum_{\ell \in \{0,1\}}
\overline{\chi_{j\ell}(h,h', k)} \mathcal{I}_{\ell}\left(\frac{1}{kz}\right) .
\end{align*}
Dividing by $\eta^6$ and applying (\ref{EtaTrans}) yields that
\begin{multline} \label{fjTrans}
f_j\left(\frac{1}{k}(h+iz)\right)=
z^{\frac{3}{2}} e^{\frac{\pi ih'}{2k}-\frac{\pi i(j+1)h}{2k} }
\sum_{\ell \in \{0,1\}}  \psi_{j\ell}(h,h', k)f_{\ell}\left(\frac{1}{k}\left(h'+\frac{i}{z}\right)
\right)  \\
+ \frac{1}{4 \sqrt{2} \pi}
z^{\frac{3}{2}}  e^{\frac{\pi ih'}{2k}-\frac{\pi i(j+1)h}{2k} }
\sum_{\ell \in \{0,1\}}  \psi_{j\ell}(h,h', k)
 \eta^{-6} \left(\frac{1}{k}
\left(h'+\frac{i}{z}\right)\right)
\mathcal{I}_{\ell}\left(\frac{1}{kz}\right)
\end{multline}
with %
\begin{align} \label{psimult}
\psi_{00}(h, h', k)&:=
\begin{cases}
- \frac{\omega_{h,k}^4\omega_{h,\frac{k}{2}}^5}{\omega_{h,\frac{k}{4}}^2}
&\text{ if } 4|k,\\[2ex]
-\frac{1}{\sqrt{2}}\frac{\omega_{h,k}^4\omega_{2h,k}^5}{\omega_{4h,k}^2}     &\text{ if
}  2\nmid k,\\
0 &\text{ if } 2\|k,
\end{cases}  \qquad
\psi_{01}(h, h', k) :=
\begin{cases}
0 &\text{ if } 4|k,\\
-\frac{1}{\sqrt{2}}   \frac{\omega_{h,k}^4\omega_{2h,k}^5}{\omega_{4h,k}^2}      &\text{
if }  2\nmid k,\\[2ex]
-\frac{\omega_{h,k}^4\omega_{h,\frac{k}{2}}^5}{\omega_{2h,\frac{k}{2}}^2}e^{\frac{\pi
    i h'}{2k}}
 &\text{ if }  2\|k,
\end{cases} \\ \nonumber
\psi_{10}(h, h', k)&:=
\begin{cases}
0                                                                                                                                                                                                                                                                                                                                 &\text{if } 4|k,\\
-\frac{1}{\sqrt{2}}\frac{\omega_{h,k}^6  \omega_{4h,k}^2}{\omega_{2h,k}}
&\text{if } 2\nmid k,\\[2ex]
- \frac{\omega_{h,k}^6
\omega_{2h,\frac{k}{2}}^2}{\omega_{h,\frac{k}{2}}^2}     &\text{if } 2\|k,
\end{cases} \qquad \,\,
\psi_{11}(h, h', k):=
\begin{cases}
-
 \frac{\omega_{h,k}^6  \omega_{h,\frac{k}{4}}^2
 }{
 \omega_{h,\frac{k}{2}}^2
 }e^{\frac{\pi ih'}{2k}}             &\text{if } 4|k,\\[2ex]
         \frac{1}{\sqrt{2}}  \frac{\omega_{h,k}^6  \omega_{4h,k}}{\omega_{2h,k}}&\text{if } 2\nmid k,\\
0                                                                                                                                                                                                                                                                                                                                         &\text{if } 2\|k.
\end{cases}
\end{align}
For later purposes we require a different representation of $I_{\ell}(x)$.
Similarly as  in \cite{BL1}, one can show
\begin{lemma}\label{nonholpart}
\noindent
We have for $x\in\C$ with Re$(x)>0$
\begin{equation}\label{thetaintegral}
\mathcal{I}_j(x)=\sum\limits_{{ g \pmod{2k}}\atop{g\equiv j\pmod{2}}}
e\left(-\frac{g^2 h'}{4k}\right)
\left(\frac{2\delta_{0,g}}{\sqrt{x}}-\frac{1}{\sqrt{2}\pi k^2 x}
\int_{-\infty}^\infty e^{-2\pi xu^2}   f_{k,g}(u)\,du\right),
\end{equation}
where $\delta_{0,g}=0$ unless $g\equiv 0 \pmod{2k}$ in which case it equals 1.\\
\end{lemma}
Note that we corrected a sign error in the statement of Lemma 4.4 of \cite{BL1}.

\subsection{Proof of Theorem  \ref{th:coefficients}}
Throughout this section, we  use the notation from Subsection \ref{TransformationSection}.
For the proof of Theorem \ref{th:coefficients}, we employ the Hardy-Ramanujan Circle Method \cite{RZ} and write for  $j\in\{0,
1\}$
\[
\widetilde{f_j}(q):=q^{\frac{j+1}{4}}f_j(\tau)=\sum_{n=0}^\infty \alpha_j(n)q^n.
\]
\noindent
By Cauchy's Theorem we have for $n>0$
\[
\alpha_j(n)=\frac{1}{2\pi i}\int_C\frac{\widetilde{f_j}(q)}{q^{n+1}}\,dq,
\]
\noindent
where $C$ is an arbitrary path inside the unit circle looping around $0$
counterclockwise. We choose the circle with radius $r=e^{\frac{-2\pi}{N^2}}$,
where we later let $N\to \infty$, and decompose it into consecutive Farey arcs of order $N$:
\[
\alpha_j(n)=\sum_{{0\leq h<k\leq N}\atop{(h, k)=1}} e^{-2\pi i n \frac{h}{k}}\int_{-\vartheta_{h, k}^{'}}^{\vartheta_{h,
k}^{''}}  \widetilde{f_j}\left(
  e^{\frac{-2\pi}{N^2}+2\pi i \frac{h}{k}+2\pi
i\phi }\right)\, e^{\frac{2\pi n}{N^2}-2 \pi in\phi}d\phi
\]
\noindent
with:
\[
\vartheta_{h, k}^{'}:=\frac{1}{k(k_1+k)},\qquad \vartheta_{h,
k}^{''}:=\frac{1}{k(k_2+k)},
\]
where $\frac{h_1}{k_1}<\frac{h}{k}<\frac{h_2}{k_2}$ are adjacent Farey fractions
in the Farey sequence of order $N$. From the theory of Farey fractions it is
known that
\[
\frac{1}{k+k_j}\leq\frac{1}{N+1}\qquad (j=1,2).
\]
\noindent
Using the transformation law (\ref{fjTrans}) and $z=k(\frac{1}{N^2}-i\phi)$, we obtain
\begin{align*}
& \alpha_j(n)=
\sum\limits_{{0\leq h<k\leq N}\atop{(h, k)=1}} e^{\frac{-2\pi ihn}{k}}
\sum_{\ell \in \{0,1\}}
\psi_{j\ell}(h,h', k) e^{\frac{\pi ih'}{2k}}
\int_{-\vartheta_{h, k}^{'}}^{\vartheta_{h, k}^{''}}
f_{\ell} \left(\frac{1}{k}\left(h'+\frac{i}{z}\right)\right)\,
e^{ \frac{2\pi z}{k}\left(n-\frac{j+1}{4} \right)} z^{\frac{3}{2}}\,d\phi\\
&
  + \frac{1}{4 \sqrt{2}\pi}
  \sum\limits_{{0\leq h<k\leq N}\atop{(h, k)=1}}e^{\frac{-2\pi
ihn}{k}}
 \sum_{\ell \in \{0,1\}}
\psi_{j\ell}(h,h', k) e^{\frac{\pi ih'}{2k}}
 \int_{-\vartheta_{h, k}^{'}}^{\vartheta_{h, k}^{''}}
\frac{\mathcal{I}_{\ell}\left(\frac{1}{k z}\right)}{\eta^6\left(\frac{1}{k}\left(h'+\frac{i}{z}\right)\right)}\,
 e^{ \frac{2\pi z}{k}\left(n-\frac{j+1}{4} \right)} z^{\frac{3}{2}}\,d\phi.
\end{align*}
We will abbreviate the first summand by $\sum_1$ and the second by  $\sum_2$.

We first consider $\sum_1$  and split of the terms with negative
exponent in the Fourier expansion as they contribute to the main term.
For this, we write
\begin{align*}
& f_0(\tau)=-\frac{1}{12}\,q^{-\frac{1}{4}}+\sum_{n>0} b_0(n) q^{n-\frac{1}{4}},\\
& f_1(\tau)=\sum_{n>0} b_1(n) q^{n-\frac{1}{2}}.
\end{align*}
We denote the contributions of the negative exponent to $\sum_1$ by $\sum_1^*$.
Using the estimates
$\vartheta_{h, k}^{'},\vartheta_{h, k}^{''}\ll\frac{1}{kN}$ and $|z|^2\ll
N^{-2}$ gives that
\[
\sum_1= \sum_1^*+O \left( N^{-\frac{5}{2}}\sum_{0\leq h<k\leq
N}\frac{1}{k}\right)
=
\sum_1^*+O\left(N^{-\frac{3}{2}}\right).
\]
To estimate $\sum_2$, we use the representation of $\mathcal{I}_{\ell}$ given in
Lemma \ref{nonholpart}.
We start with the contribution of the first summand in the representation of $\mathcal{I}_{\ell}$ (only occuring for
$\ell=0$ and  $g \equiv 0 \pmod{2k}$). Splitting of the non-principal terms yields as before an error of $N^{-\frac32}$.
 To estimate the remaining terms of $\sum_2$, we
  aim to estimate integrals of the shape
\[
\mathcal{I}_{k, g, b}(z):=e^{\frac{2\pi b}{kz}}z^{\frac{5}{2}}
\int_{-\infty}^\infty e^{-\frac{2\pi u^2}{kz}}f_{k, g}(u)\,du.
\]
Similarly to the case of Fourier expansions we are interested in the ''principal integral part''
contribution. To be more precise, we let for  $b>0$ and $g\in\Z$,
\[
\mathcal{J}_{k, g, b}(z):=e^{\frac{2\pi b}{kz}}z^{\frac{5}{2}}
\int_{-\sqrt{b}}^{\sqrt{b}}e^{-\frac{2\pi u^2}{kz}}f_{k, g}(u)\, du.
\]
Similarly as in \cite{BM}, we may show:
\begin{lemma}\label{nonholestimate}
As $z\to\infty$ we have for $-k<g \leq k$:
\begin{itemize}
\item[(1)]
If $b\leq 0$, then
\begin{equation} \label{IEstimate}
\left|\mathcal{I}_{k, g, b}(z)\right|\ll
|z|^{\frac{5}{2}} \times
\left\{
\begin{array}{ll}
\frac{k^2}{g^2} & \text{if } g\not=0,\\
1&\text{if } g=0.
\end{array}
\right.
\end{equation}
\item[(2)]
If $b>0$, then
\[
\mathcal{I}_{k, g, b}(z)=\mathcal{J}_{k, g, b}(z)+\mathcal{E}_{k, g, b}(z),
\]
where the error $\mathcal{E}_{k,g,b}$ satisfies the same estimate as $\mathcal{I}_{k, g, b}$ in
(\ref{IEstimate}).
\end{itemize}
\end{lemma}

\begin{proof}
Recall that $Re\left(\frac{1}{z}\right)\geq\frac{k}{2}$.
We use the estimate
\[
\Bigg|\sinh\left(\frac{\pi u}{k}-\frac{\pi ig}{2k}\right)\Bigg|=\Bigg|\cosh\left(\frac{\pi u}{k}-\pi i\left(\frac{g}{2k}+\frac{1}{2}\right)\right)\Bigg|
\geq \Big|\sin\left(\frac{\pi g}{2k}\right)\Big|.
\]
Now for $0\leq x \leq\frac{\pi}{2}$, the function $\frac{\sin(x)}{x}$ is bounded from below, thus for $-k<g\leq k$, $g\neq 0$
\[
\Bigg|\frac{1}{\sin\left(\frac{\pi g}{2k}\right)^2}\Bigg|\ll \frac{k^2}{g^2}.
\]
Moreover
\[
\Bigg|\frac{1}{\sinh^2(x)}-\frac{1}{x^2}\Bigg|=\frac{1}{x^2}-\frac{1}{\sinh^2(x)}\leq 1.
\]
We now define $h_{g, k}$ as
\[
h_{g, k}:=
\begin{cases}
\frac{k^2}{g^2}  &\quad\text{if } -k<g\leq k, g\neq 0,\\
1                &\quad\text{if } g=0.
\end{cases}
\]
Then by the above
\begin{equation*}\label{fkestimate}
|f_{k, g}(u)|\leq h_{g, k}.
\end{equation*}
We now first assume that $b\leq 0$. Then
\[
|I_{k, g, b}(z)|\leq |z|^{\frac{5}{2}}h_{g, k} \int_{-\infty}^\infty e^{-\frac{2\pi u^2}{k}Re\left(\frac{1}{z}\right)}du
\ll |z|^{\frac{5}{2}}h_{g, k}\sqrt{\frac{k}{Re\left(\frac{1}{z}\right)}}\ll |z|^{\frac{5}{2}}h_{g, k}
\]
which gives the claim for $b\leq 0$. The case $b>0$ works similarly.
\end{proof}

The contribution of non-principal part of the remaining terms of $\sum_3$ may be estimated against
a constant times
$$
 \sum_{h,k} \frac{1}{k}   \int_{-\vartheta_{h, k}^{'}}^{\vartheta_{h, k}^{''}}
\sum_{g} \left|  \mathcal{I}_{k,g,0}(z)\right| d \phi
\ll \sum_{h,k}  \frac{1}{k}   \int_{-\vartheta_{h, k}^{'}}^{\vartheta_{h, k}^{''}}|z|^{\frac52}
\left( 1 + \sum_{g=1}^{k} \frac{k^2}{g^2}
\right) d \phi \ll N^{-\frac32}.
$$
 In the terms coming from the principal part, we may similarly truncate the
integral to lead $\mathcal{J}_{k,g,\frac14}(z)$.
Combining the above, we have shown
that
\[
\alpha_j (n)=S_1+S_2+S_3+O\left(N^{-\frac{3}{2}}\right)
\]
with
\begin{eqnarray*}
S_1&:=&-\frac{1}{12} \sum\limits_{{0\leq h<k\leq N}\atop{(h, k)=1}}
 e^{\frac{-2\pi ihn}{k}}
\psi_{j0}
(h, h', k)
\int_{-\vartheta_{h, k}'}^{\vartheta_{h,
k}''}
e^{\frac{2\pi z}{k}\left(n-\frac{(j+1)}{4}\right)+\frac{\pi}{2kz}} z^{\frac{3}{2}}d\phi,
\\
S_2&:=&
\frac{1}{2\sqrt{2} \pi} \sum\limits_{{0\leq h<k\leq N}\atop{(h, k)=1}}
\sqrt{k}  e^{\frac{-2\pi
ihn}{k}}
\psi_{j0}
(h, h', k) \int_{-\vartheta_{h, k}'}^{\vartheta_{h, k}''}
e^{\frac{2\pi z}{k}\left(n-\frac{j+1}{4}\right)+\frac{\pi}{2kz}} z^2 d\phi,
\\
S_3&:=&-\frac{1}{8\pi^2}\sum\limits_{{0\leq h<k\leq N}\atop{(h,
k)=1}}
 \frac{1}{k}\
 e^{\frac{-2\pi
ihn}{k}} \hspace{-3ex}
\sum_{\substack{\ell \in \{0,1\}\\ \substack{ -k< g\leq k\\ g\equiv \ell \pmod{2}}}}
\hspace{-3ex}
\psi_{j\ell}
(h, h', k)
  e\left(\frac{-g^2
h'}{4k}\right)
\int_{-\vartheta_{h, k}'}^{\vartheta_{h, k}''}e^{\frac{2\pi
z}{k}\left(n-\frac{j+1}{4}\right)} \mathcal{J}_{k, g, \frac{1}{4}}(z)d\phi
   .
\end{eqnarray*}
We next write the path of integration in a symmetrized way
\[
\int_{-\vartheta_{h, k}^{'}}^{\vartheta_{h, k}^{''}}
=\int_{-\frac{1}{kN}}^{\frac{1}{kN}}
- \int_{-\frac{1}{kN}}^{-\frac{1}{k(k+k_1)}}
- \int_{\frac{1}{k(k+k_2)}}^{\frac{1}{kN}}\ .
\]
The second and third term contribute to the error term and may be estimated as
before. To finish the proof, we require estimates for integrals of the form $(r>0)$
\[
\mathcal{I}_{k, r, n, m}:=\int_{-\frac{1}{kN}}^{\frac{1}{kN}}z^r
e^{\frac{2\pi}{k}\left(nz+\frac{m}{z}\right)}d\phi.
\]
In a standard way  (we refer the reader to \cite{B2}  for the details) one may show that
\[
\mathcal{I}_{k, r, n,
m}=\frac{2\pi}{k}\left(\frac{m}{n}\right)^{\frac{r+1}{2}}I_{r+1}\left(\frac{4\pi}{k}\sqrt{nm}\right)
+O\left(\frac{1}{kN^{r+1}}\right)\,.
\]
Inserting this bound into the $S_i$ and letting $N \to \infty$ now
easily gives the claim. $\square$

\providecommand{\href}[2]{#2}\begingroup\raggedright

\end{document}